\documentclass[12pt]{article}
\usepackage{amssymb}
\usepackage{latexsym}
\usepackage[all]{xy}
\newcounter{mytab}
\newenvironment{mytab}[1]{
\refstepcounter{mytab}
\begin{center}
Table \arabic{mytab}. #1
\end{center}
}{}
\newtheorem{THM}{Theorem}[section]

\newtheorem{LEM}{Lemma}[section]

\newcommand{\oao}{o-o }
\newcommand{\oae}{o-e }
\newcommand{\eae}{e-e }
\begin{document}
%
%

\title{Locally Toroidal Polytopes and Modular Linear Groups}
 
\author{B. Monson\thanks{Supported by NSERC of Canada Grant \# 4818}\\
University of New Brunswick\\
Fredericton, New Brunswick, Canada E3B 5A3
\and and\\[.05in]
Egon Schulte\thanks{Supported by NSA-grants H98230-05-1-0027 and
H98230-07-1-0005}\\
Northeastern University\\
Boston, Massachussetts,  USA, 02115}

\date{ \today }
\maketitle

\begin{abstract}
When the standard  
representation of a crystallographic Coxeter group $G$ 
(with string diagram) is reduced modulo the integer $d \geq 2$,
one obtains a finite group $G^d$ which is often the 
automorphism group of an abstract regular polytope. 
Building on earlier work
in the case that $d$ is an odd prime, we here                                                                                                                                                                    develop  
methods to handle composite moduli and completely describe the 
corresponding modular polytopes when
$G$ is of spherical or Euclidean type.
Using a modular variant of the quotient criterion, we then describe 
the locally toroidal
polytopes provided by our construction, most of which are new.



\end{abstract}

\section{Introduction}
\label{intro}

Our fascination with the regular polytopes is due not only to their 
visual appeal and charm, but also to the fact that their symmetry
groups appear in such varied and unexpected places. 
In a recent series of papers, for example, the authors established
the basic machinery needed to describe a large class of polytopes
whose automorphism groups typically have small index
in some finite orthogonal group  
(see \cite{monsch1,monsch2,monsch3}). 
Indeed, in our analysis there  we  often had to exploit quite subtle
properties of the orthogonal group $O(n,p,\epsilon)$
on an $n$-dimensional vector space over $\mathbb{Z}_p$, 
where $p$ is an odd prime.
Here we take a bit of a detour 
and consider instead the possibilities released by more generally
working over the ring $\mathbb{Z}_d$, with any modulus $d \geq 2$.
(The rank $4$ polytopes described in \cite{monwei1,monwei2} involve 
an analogous excursion into the domains of Gaussian and Eisenstein
integers; and, of course, the related  idea of  constructing the automorphism
group of a regular map by modular reduction is natural and
well established; see \cite{wilson}, for example.)

Our main goal is to extend  previous results on 
\textit{locally toroidal polytopes},  as provided by our
construction  
\cite[\S 4]{monsch3}. 
To that end, in Sections~\ref{absgps} and \ref{crysgps} we describe 
the modular reduction of a crystallographic Coxeter group $G$ 
with string diagram. In Sections ~\ref{modsph} and \ref{euctype}
we completely describe what happens when $G$ is of spherical or Euclidean type.
Finally, after proving a useful quotient criterion (Theorem~\ref{interquotA}), 
we discuss in Section~\ref{loctor} various new families of locally toroidal
polytopes, mainly in ranks $5$ and $6$.

\section{Abstract regular polytopes and Coxeter groups}
\label{absgps} 
 
Let us begin  with a brief review of some  key properties of
abstract regular polytopes,
referring to \cite{arp} for details. 
An {\em (abstract) $n$-polytope} $\mathcal{P}$ is a partially ordered set
with a strictly monotone rank function having range $\{-
1,0,\ldots,n\}$. An element $F \in \mathcal{P}$ with ${\rm rank}(F)=j$ is
called a $j$-{\it face}; typically $F_j$ will indicate a $j$-face;
$\mathcal{P}$ has a unique least face
$F_{-1}$ and   unique greatest face $F_n$. 
Each maximal
chain or {\it flag} in $\mathcal{P}$ must contain $n+2$ faces. Next, 
$\mathcal{P}$ must satisfy a 
\textit{homogeneity property}\,: whenever $F < G$ with ${\rm rank}(F)=j-1$ and
${\rm rank}(G)=j+1$, there are exactly two $j$-faces $H$ with
$F<H<G$, just as happens for convex $n$-polytopes. 
It follows that for $0 \leq j \leq n-1$ and any flag $\Phi$, 
there exists a unique \textit{adjacent} 
flag ${}^j\Phi$, differing from $\Phi$ in just the
rank $j$ face. With this notion of adjacency the flags of $\mathcal{P}$
form a \textit{flag graph}.
The final defining property of $\mathcal{P}$ is that 
the flag graph for each
section must be connected, so that $\mathcal{P}$ is 
\textit{strongly flag--connected}. Recall here that whenever 
$F \leq G$ are faces of ranks $j \leq k$
in $\mathcal{P}$,  then the \textit{section} of $\mathcal{P}$ determined by $F$ and $G$ is given by  
$G/F := \{ H \in \mathcal{P}\, | \, F \leq H \leq G \}$. In fact, this
is a ($k-j-1$)-polytope in its own right.

Naturally, the symmetry of $\mathcal{P}$ is exhibited by its
\textit{automorphism group} $\Gamma(\mathcal{P})$, containing all
order preserving bijections on $\mathcal{P}$. Henceforth, we shall 
consider only \textit{regular polytopes} $\mathcal{P}$, for which  
$\Gamma(\mathcal{P})$ is, 
by definition, transitive on flags. 
Clearly a regular $n$-polytope $\mathcal{P}$ must have
all sorts of local combinatorial symmetry.
In particular,   $\mathcal{P}$ will be  \textit{equivelar}
of some type $\{p_1, \ldots , p_{n-1}\}$, where
$2 \leq p_j \leq \infty$;  this means that for 
each fixed $  j \in \{1, \ldots,n-1\}$ and  
each pair of incident faces $F$ and $G$  in 
$\mathcal{P}$, with
$\mbox{rank} (F) = j-2$ and $\mbox{rank} (G) = j+1$,  the rank 2
section $G/F$ has the structure of a $p_j$-gon (independent of choice
of $F < G$).  Thus, each 2-face (polygon) of $\mathcal{P}$ is isomorphic to
a $p_1$-gon, and  in every 3-face of $\mathcal{P}$, each
0-face is surrounded by an alternating cycle of $p_2$ edges and $p_2$
polygons, etc. 

To further understand the structure of $\Gamma(\mathcal{P})$ when
$\mathcal{P}$ is regular, we  
fix a
\textit{base flag} $\Phi = \{F_{-1},F_0,\ldots,F_{n-1},F_n\}$, with
$\mbox{rank}\; (F_j) = j$.  For $0 \leq j \leq n-1$,
let $\rho_j$ be the (unique) automorphism with
$\rho_j (\Phi) = {}^j\Phi$.  If $\mathcal{P}$ is regular, then
$\Gamma({\cal P})$ is generated by 
$\rho_0,\rho_1,\ldots,\rho_{n-1}$,
which are involutions satisfying 
 at least the relations 

\begin{equation}
\label{genrels1}
\rho_j^2 = (\rho_{j-1} \rho_j)^{p_{j}} = (\rho_{i} \rho_j)^{2} 
= 1,\;\;\;\;0 \leq i,j \leq n-1,
\; |j-i| \geq 2\; 
\end{equation}
(with $j \geq 1$ for $\rho_{j-1} \rho_j$). Also, an {\em intersection
condition} on standard subgroups holds: 
\begin{equation}
\label{geninter}
\langle \rho_i\,|\, i \in I \rangle \cap \langle \rho_i\,|\,i \in J \rangle =
\langle \rho_i\,|\, i \in I \cap J \rangle
\end{equation}
for all $I,J \subseteq \{0,\ldots,n-1\}$. 
In short, $\Gamma({\cal P})$
is a very particular quotient of the  Coxeter group
$G = [p_1, \ldots, p_{n-1}]$,  
whose diagram  is a string  with branches labelled
$p_1, \ldots , p_{n-1}$. (We allow $p_j=2$, in which case the `string' is
disconnected.)
Conversely, given any group $\Gamma = \langle \rho_0,\ldots,\rho_{n-1}
\rangle$ generated by involutions and
satisfying (\ref{genrels1})  and (\ref{geninter}), 
one may construct a
polytope ${\cal P}$ with $\Gamma({\cal P}) = \Gamma$ (see
\cite[Theorem 2E11]{arp}).
We then  say that  $\Gamma({\cal P})$ is a {\em string C-group}.
Since $\mathcal{P}$ can be uniquely
reconstructed from  $\Gamma(\mathcal{P})$,  we may therefore shift our focus
to an appropriate  class of groups of particular interest.

Recall that if $\mathcal{P}_1 $ and $\mathcal{P}_2$ are regular $n$-polytopes with $n\geq 2$, then 
$\langle\, \mathcal{P}_1\, , \,\mathcal{P}_2 \, \rangle$ denotes the class of all
regular $(n+1)$-polytopes whose facets are isomorphic to 
$\mathcal{P}_1$ and whose vertex-figures are isomorphic to
$\mathcal{P}_2$. If this class is non-empty, then it contains a 
\textit{universal} regular $(n+1)$-polytope, denoted 
$\{ \, \mathcal{P}_1\, , \,\mathcal{P}_2 \, \}$,
which covers any other polytope in the class \cite[Th. 4A2]{arp}.

Let us look more closely at  the abstract Coxeter group
$G = [p_1, \ldots, p_{n-1}]$, which   is itself a string C-group
with respect to the usual generators and which may well be infinite. 
The corresponding polytope
$\{ p_1, \ldots ,p_{n-1}\}:= \mathcal{P}(G)$ is universal in a more local sense,
as described in \cite[Th. 3D5]{arp}. 

Now, like any finitely generated Coxeter group, 
$G$ can be identified with its image under the standard 
faithful representation in real $n$-space $V$ \cite[Cor. 5.4]{humph}.
Consequently, we may suppose 
$G =  \langle r_0, \ldots , r_{n-1} \rangle$ to be the
\textit{linear Coxeter group} generated by certain
reflections $r_j$ on $V$. In fact, these
reflections leave invariant a symmetric bilinear form $x \cdot y$ on $V$,
so that $G$ is a subgroup of the corresponding orthogonal group
$O(V) \subset GL(V)$. (Note, however,  that $x \cdot y$ is positive definite 
if and only if
$G$ is finite \cite[Th. 6.4]{humph}.)
We shall let $e$ denote the identity in the group $GL(V)$.

Recalling our earlier description of the regular $n$-polytope $\mathcal{P}$,
we now have an epimorphism
\begin{eqnarray*} 
G & \rightarrow & \Gamma(\mathcal{P})\\
r_j & \mapsto &\rho_j \;.
\end{eqnarray*}
Intuitively then, we may think of regular polytopes as having maximal
reflection symmetry.

\section{Crystallographic Coxeter groups and their modular reductions}
\label{crysgps}

Now let us specialize.  We   say that the linear Coxeter group
$G$ is \textit{crystallographic} (with
respect to the standard representation) if it leaves invariant
some \textit{lattice} $\sum_{j=0}^{n-1} \mathbb{Z} b_j$ 
generated by a basis $\beta = \{b_j\}$ 
for $V$.  As described
in  \cite{max2}  or \cite[Prop. 4.1]{monsch1}, there is no loss of 
generality in assuming that $\beta$ is a \textit{basic system} for 
$G$, meaning that each  $b_j$ is a \textit{root} for the corresponding 
reflection $r_j$. Thus,  

\begin{equation}\label{refln}
r_i(b_j) = b_j + m_{i\, j}b_i
\end{equation}
for certain \textit{Cartan integers} $m_{i \, j}$, $0 \leq i,j \leq n-1$,
with   all
$m_{i\,i} = -2$ and $m_{i\, j} = 0$ for $|i -j| \geq 2$.

Now recall that the string Coxeter group
$G = [p_1, \ldots, p_{n-1}]$ is crystallographic if and only if 
all $p_j \in \{2,3,4,6,\infty\}$ \cite[Prop. 4.1(c)]{monsch1}. 
If the corresponding   \textit{Coxeter diagram} $\Delta_{c}(G)$
is connected, then $G$ admits only a finite number of essentially
distinct basic systems $\beta$. 
As we observed in \cite[\S 4]{monsch1},
each basic system and   corresponding lattice
can be encoded in a new diagram $\Delta(G)$,
a variant of   $\Delta_{c}(G)$. Briefly, the 
branches of  $\Delta(G)$ are no longer labelled; instead, 
each node $j$ of $\Delta(G)$
is  labelled  by the real number $b_j^2 = b_j \cdot b_j$. 
Each subdiagram  on two nodes $i$ and $j$
must then be one of those appearing 
in Table~\ref{BasicSystems} below.

\begin{center}
\begin{tabular}{c|c|c}\label{dihedgps}
Period of $r_i r_j$ & Subdiagram  on nodes & Cartan integers \\
& $i$ (left), $j$ (right) & $m_{i j} , m_{j i}$ \\\hline\hline
 $2$& $ \stackrel{a}{\bullet}\!\;\;\;\;\;\;\!
	\stackrel{c}{\bullet}$&  $0, 0$\\\hline
 $3$ & $\stackrel{a}{\bullet}\!\frac{}{\;\;\;\;\;\;}\!
	\stackrel{a}{\bullet}$& $1, 1$ \\\hline
 $4$& $\stackrel{a}{\bullet}\!\frac{}{\;\;\;\;\;\;}\!
	\stackrel{2a}{\bullet}$& $2, 1$\\\hline
 $6$& $\stackrel{a}{\bullet}\!\frac{}{\;\;\;\;\;\;}\!
	\stackrel{3a}{\bullet}$& $3, 1$\\\hline
 $\infty$ & $\stackrel{a}{\bullet}\!\frac{}{\;\;\;\;\;\;}\!
	\stackrel{4a}{\bullet}$& $4, 1$ \\ \hline
 $\infty$ & $\stackrel{a}{\bullet} =\!\!=\!\!= \stackrel{a}{\bullet}$
 	& $2,2$ \\ 
\end{tabular}

\vspace*{3mm}
 
\end{center}
\begin{mytab}{Possible diagrams for dihedral subgroups
$\langle r_i , r_j \rangle$ of $G$}\label{BasicSystems}\end{mytab}
For each $i \neq j$, we have $m_{i j} m_{j i} = 4\cos^2(\pi/p_{i j})$,
where $p _{i j}$ is the period of the \textit{rotation} $r_i r_j$. 
(In particular, $p_{j-1, j} = p_j$.) Note that
nodes $i$ and $j$ must be clearly distinguished, say as left and right in the 
Table~\ref{BasicSystems}, whenever $m_{i j}  \neq m_{j i}$.
By suitably rescaling the node labels on each connected 
component of $\Delta(G)$, we can assume that these labels are a 
set of relatively prime positive integers. As a familiar example, consider the
usual tessellation $\mathcal{P}$ of the Euclidean plane by congruent squares. 
Then $\mathcal{P}$ is an infinite  regular $3$-polytope,
and $G = [4,4] \simeq \Gamma(\mathcal{P})$ admits the diagrams
\begin{equation}\label{44tors}
\stackrel{1}{\bullet}\!\frac{}{\;\;\;\;\;\;}\!
	\stackrel{2}{\bullet}\!\frac{}{\;\;\;\;\;\;}\!
	\stackrel{1}{\bullet}\; , \;\;
	\stackrel{1}{\bullet}\!\frac{}{\;\;\;\;\;\;}\!
	\stackrel{2}{\bullet}\!\frac{}{\;\;\;\;\;\;}\!
	\stackrel{4}{\bullet}\;  \mbox{   and   }\;
	\stackrel{2}{\bullet}\!\frac{}{\;\;\;\;\;\;}\!
	\stackrel{1}{\bullet}\!\frac{}{\;\;\;\;\;\;}\!
	\stackrel{2}{\bullet}
\;. 
\end{equation}

Having fixed such a basic system for a crystallographic
Coxeter group $G = [p_1, \ldots , p_{n-1}]$, we can reduce 
$G$ modulo any integer $s \geq 2$: the natural 
epimorphism $\mathbb{Z} \rightarrow \mathbb{Z}_s$ induces a homomorphism
of $G$ onto a  subgroup $G^s$ of $GL_n(\mathbb{Z}_s)$, 
the group of $n \times n$ invertible matrices over  $\mathbb{Z}_s$.
Our hope, of course, is that the finite group $G^s$ will be the 
automorphism group of a regular $n$-polytope.
(In \cite{monsch1,monsch2, monsch3}   we examined such groups in the case 
that $s$ is an odd prime, so as to 
exploit the structure of orthogonal groups over 
finite fields.)

We shall often abuse notation by referring to the modular images of objects by 
the same name (such as $r_i$,   $e$, $b_i$, $V$, etc.). In particular, 
$\{b_i\}$ will denote the standard basis for $V = \mathbb{Z}_s^n$, which in general
we must now view as a free module over the ring $\mathbb{Z}_s$.   
We shall see in Lemma~\ref{redmodgen} that 
$r_{i}$  usually continues to act as a reflection after reduction; 
 in any case, we can  
compute it using (\ref{refln}). However, the situation for metrical quantities
such as $b_i \cdot b_j$, a rational number which occasionally
has denominator $2$, is more intricate \cite[Eq. 10]{monsch1}. Nevertheless,
at least when $\mbox{gcd}(6,s) =1$, we can interpret
$G^s$ as a subgroup of the orthogonal group $O( \mathbb{Z}_s^n)$ for 
the  symmetric bilinear form   
defined on $\mathbb{Z}_s^n$ by means of the Gram matrix $[b_i \cdot b_j]$.
Moreover, we can then write 
$$ r_i(x) = x - 2 \, \frac{x \cdot b_i}{b_i \cdot b_i}\, b_i$$
since $b_i^2$ will be invertible $\pmod{s}$.
In our earlier work with prime moduli, these issues were a concern 
only for `non-generic'
groups, where  $s = 3$ and $G$ has some period $p_j = 6$. 
Here, with more  general moduli, the analysis is   more 
complicated. It often happens, for instance,  that
the group $G^s$   depends essentially on the choice of 
 basic system and  the corresponding diagram 
$\Delta(G)$.
For example, for the modulus $s=4$, the group $G^4$ corresponding to the 
three diagrams in (\ref{44tors}) has order $32$, $128$ and $64$, respectively.
These are, in fact, the automorphism groups of the regular toroidal maps
$\{4,4\}_{(2,0)}$, $\{4,4\}_{(4,0)}$  and $\{4,4\}_{(2,2)}$
(see Table~\ref{Cubics} below). 

Clearly, we must now confront a crucial question: 
when is $G^s = \langle r_0, \ldots , r_{n-1} \rangle^s$ a string $C$-group 
(i.e. the automorphism group of a finite, 
abstract regular $n$-polytope $\mathcal{P} = \mathcal{P}(G^s)$)?
Unfortunately,  we cannot provide anything like a comprehensive
answer here. 
Instead, for classes of groups $G$ of particular interest,
we shall have to rely more on 
ad hoc techniques 
than we did for  prime moduli, without trying to
exploit in any deep way the structure of orthogonal groups over  
general rings. Occasionally, we employ
GAP \cite{gap} to settle `small' cases.

Certainly, the generators $r_j$ of $G^s$   \textit{satisfy}
the Coxeter-type relations  inherited from $G$ (see 
(\ref{genrels1}), with $\rho_j$ replaced by $r_j$). However,  
before confronting the intersection condition (\ref{geninter}) for 
$G^s$, we must take a closer look.
For example, it might happen that 
$r_j =e \pmod{s}$.

\noindent\textbf{Notation}. We say that node $i$ of $\Delta(G)$ is 
\textit{\eae} if both Cartan integers $m_{i, i-1}$ and $m_{i, i+1}$ are   even; 
\textit{\oae} if just one of the integers is even; and \textit{\oao} if 
both are odd.
For  the terminal nodes $0$ and $n-1$ on the string we shall
agree that $m_{0 , -1}=m_{n-1 , n} = 0$.

Note that end nodes can never be \oao.
Likewise,  a node is   \eae if 
it is 
labelled $a$, while any adjacent nodes 
are labelled
$4a$, $2a$ or $a$ (after a double branch), as in 

$$ \ldots \frac{}{\;\;}\!\!\stackrel{2a}{\bullet}\!\frac{}{\;\;\;\;\;\;}\!
	\stackrel{a}{\bullet}\!\frac{}{\;\;\;\;\;\;}\!
        \stackrel{2a}{\bullet}\!\!\frac{}{\;\;} \ldots \;\;
,\;\;
\stackrel{a}{\bullet} =\!\!=\!\!= 
        \stackrel{a}{\bullet}\!\!\frac{}{\;\;} \ldots
\;\;,\;\;\ldots
\frac{}{\;\;}\!\!\stackrel{2a}{\bullet}\!\frac{}{\;\;\;\;\;\;}\!
        \stackrel{a}{\bullet}\, ,
\;\;\mbox{\rm etc.}\;$$
Typical \oae nodes are the middle nodes in the  subdiagrams
$$
\ldots \frac{}{\;\;}\!\!\stackrel{3a}{\bullet}\!\frac{}{\;\;\;\;\;\;}\!
        \stackrel{a}{\bullet}\!\frac{}{\;\;\;\;\;\;}\stackrel{2a}{\bullet}\!\!\frac{}{\;\;}\!\ldots
\;\;\mbox{\rm or}\;\;\ldots
\frac{}{\;\;}\!\!\stackrel{a}{\bullet} =\!\!=\!\!= 
        \stackrel{a}{\bullet}\!\frac{}{\;\;\;\;\;\;}\stackrel{c}{\bullet}\!\!\frac{}{\;\;}\!\ldots$$
(where the  integer label $c$ divides $a$).
Let us  now summarize basic  properties of the generators
$r_i$ for $G^s$. 
Using  (\ref{refln}), the calculations are straightforward, if a bit involved.

\begin{LEM}\label{redmodgen}
Let $G  = \langle r_0, \ldots , r_{n-1} \rangle \simeq [p_1, \ldots , p_{n-1}]$
be any crystallographic linear Coxeter group with 
string diagram. Suppose $s \geq 2$,
and reduce $G$ modulo $s$. Then

{\rm (a)}\,\,Each $r_i \in G^s$ has period $2$, except that $r_i = e$
when $s=2$ and  node $i$ of $\Delta(G)$ is \eae.

{\rm (b)}\, $r_i$ and $r_j$ commute in $G^s$ when $i < j-1$.

{\rm (c)}\, Suppose $p_i= 2, 3, 4$ or $6$. If $s>2$, then 
$r_{i-1} r_i$ has period $p_i$ in $G^s$
(unchanged from characteristic $0$). 

Now let $s=2$.  If  $p_i = 3$ or $6$, 
the  period of $r_{i-1} r_i$ is always $3$. If  $p_i = 4$, the
period   collapses to $2$ if and only if one  of nodes $i-1$ or $i$ is \eae. 
For $p_i = 2$, the period collapses to $1$ if and only if
both nodes are \eae (so that $r_{i-1} = r_i = e$). 

{\rm (d)}\, Suppose $p_i=\infty$. Then $r_{i-1} r_i$
has period $s$ in $G^s$, except in the following
cases, each when $s$ is even: for the subdiagram 
$\frac{}{\;\;}\!\!\stackrel{a}{\bullet} =\!\!=\!\!= 
\stackrel{a}{\bullet}\!\!\frac{}{\;\;}$, the period becomes
$\frac{s}{2}$ when both nodes\\ are  \eae;  
for the subdiagram 
$\frac{}{\;\;}\!\!\stackrel{a}{\bullet}
\!\frac{}{\;\;\;\;\;\;}\!\stackrel{4a}{\bullet}\!\!\frac{}{\;\;}$ the
period becomes $2s$ when the node labelled $a$ is \oae.
\end{LEM}

\noindent\textbf{Remarks}. In the typical case, when 
 all $r_i$ have period $2$, we say that $G^s$
is a \textit{string group generated by involutions}. Even for modulus $s=2$, it is quite
possible that all $r_i$ be involutions
(though not geometrical reflections), so long as $\Delta(G)$ 
has special features, as explained later.
Assuming now that all $r_i$ are involutions, we conclude that
$G^s$ is a string C-group  if and only if it
satisfies  the   intersection condition (\ref{geninter}), with $r_i = \rho_i$.
Our main problem is therefore to determine when $G^s$ satisfies
(\ref{geninter}). 

We hinted earlier at the definite advantages of working
with prime moduli.  For a composite modulus $s$,  we would at least hope that $G^s$ somehow splits according to the prime decomposition of $s$. 
However, our hopes for a simple approach are dashed by examples 
such as the following. 
Let $G \simeq [4,6,4]$ be the group with diagram

$$\stackrel{2}{\bullet}\!\frac{}{\;\;\;\;\;\;}\!
\stackrel{1}{\bullet}\!\frac{}{\;\;\;\;\;\;}\!
\stackrel{3}{\bullet}\!\frac{}{\;\;\;\;\;\;}\!
\stackrel{6}{\bullet}
\;.$$
First of all, we find for $p=2$ that $G^2$ is a string $C$-group of order
$96$. The middle rotation order collapses and we actually obtain the group
for the universal locally projective polytope
$\{\, \{4,3\}_3\, , \, \{3,4\}_3\, \}$.
For $p=3$ we get a group $G^3$
of order $5184$ for a self-dual polytope of type
$\{4,6,4\}$ (see \cite[Eq. (33)]{monsch2}).

Now for modulus $s=6$ we find that $G^6$ has order
$248832 = \frac{1}{2} (96 \times 5184)$. But  the intersection
condition fails, since
$\langle r_1 , r_2 \rangle^6$ has index $3$ in
$\langle r_0, r_1 , r_2 \rangle^6 \cap \langle r_1 , r_2, r_3 \rangle^6$.
In other words, the polytopality of $G^s$ is not determined through the
prime factorization of $s$. Since, in the end, we are more concerned with 
locally toroidal groups $G$,  which do fall to a  more direct  attack, we
shall mainly ignore the prime factorization of $s$.  (We note, however, that 
precisely that approach worked in \cite{monwei1,monwei2}. But for the 
$4$-polytopes considered there, the rotation groups were covered by special linear groups over
certain rings of algebraic integers; and the resulting modular groups do split
according to the prime factorization.) 

Before proceeding, let us set down  some useful notation.
For any $J \subseteq \{0, \ldots, n-1\}$, we let 
$G_J^s := \langle r_j \,|\, j \not\in J \rangle $;  in particular, for 
$k,l \in \{0, \ldots n-1\}$ we let 
$G_k^s := \langle r_j \,|\, j \neq k\rangle$ and 
$G_{k,l}^s := \langle r_j \,|\, j \neq k,l \rangle$. 
 We also let $V_J$ be the submodule of 
$V = \mathbb{Z}_s^n$ spanned by $\{b_j \!\mid\! j \not\in J \}$, and 
similarly for $V_k, V_{k,l}$. Note that $V_J$ is $G_J^s$-invariant. In particular, 
$G_j^s$ acts on $V_j$, for $j=0$ or $n-1$; however, this action need not 
be faithful
(see   \cite[Lemma 3.1]{monsch2}).

\section{Modular polytopes of spherical type}
\label{modsph}

When $G = [p_1, \ldots , p_{n-1}]$ is finite, the invariant form
$x \cdot y$  on real $n$-space $V$ 
is positive definite, so that $G$ acts 
in a natural way on any sphere $\mathbb{S}^{n-1}$ with centre $o \in V$.
Accordingly, we also say that $G$ is of \textit{spherical type}.
If the spherical group $G$  has a connected diagram, 
then up to isomorphism $\mathcal{P}(G)$
is one of the familiar convex regular $n$-polytopes \cite[\S 5-6]{monsch1}. 
After central projection, such polytopes  
can usefully  be viewed as regular spherical tessellations of the 
circumsphere $\mathbb{S}^{n-1}$.

In \cite[\S 5-6]{monsch1} we showed that $G \simeq G^p$, for any odd 
prime modulus $p$ and crystallographic string Coxeter group $G$ of spherical type 
and in any rank $n \geq 1$.  When $s$ is divisible by an odd prime $p$, the 
natural epimorphisms
$$G \rightarrow G^s \rightarrow G^p$$
immediately give $G \simeq G^s$. Here we take a 
different approach,
working explicitly with the underlying representation of the spherical 
group $G$ in 
$GL_n(\mathbb{Z})$. We confirm  that $G^s \simeq G$ for \textit{any} 
modulus $s \geq 3$,
and sometimes even for $s=2$. 
(For $n =1, 2$  such isomorphisms  follow at once from
Lemma~\ref{redmodgen}.)  
However, since the actual calculations  for general rank $n$ are quite 
tiresome, we shall simply summarize the 
results, with brief comments, for each of 
the relevant families of spherical groups. In fact, to serve later 
applications, we must generalize a little and consider how a spherical group
can  embed as a string subgroup of  
some group $G$ of higher rank. In other words, we consider  certain
\textit{spherical subdiagrams} of $\Delta(G)$.

 Of course, when
$G^s \simeq G$ we also know the structure of the modular polytope
$\mathcal{P}(G^s)$, which is merely a copy of $\mathcal{P}(G)$.

\noindent 
\textbf{(a) The group of the $m$-simplex: $A_m \simeq S_{m+1}$, for $ m \geq 1$}. 

Here, for some label $a \geq 1$, $\Delta(G)$ has  the subdiagram
\begin{equation}\label{Amsys} 
\ldots \frac{}{\;\;\;}\!
\stackrel{a}{\bullet}\!
\frac{}{\;\;\;\;\;\;}\!
\stackrel{a}{\bullet}\!
\frac{}{\;\;\;}\;  \ldots \; \frac{}{\;\;\;}\;
\stackrel{a}{\bullet}\!
\frac{}{\;\;\;\;\;\;}\!
\stackrel{a}{\bullet}\!\frac{}{\;\;\;}\;\ldots  , 
\end{equation}
on $m$ consecutive nodes $j, \ldots, j+m-1$. For all $m \geq 2$ and all $s \geq 2$,
 we then have
$$ \langle r_j, \ldots , r_{j+m-1} \rangle^s \simeq A_m\;.$$
Part (c) of Lemma~\ref{redmodgen} provides the base step of an induction on $m \geq 2$.
As in \cite[\S 6.1]{monsch1}, we then exploit the contragredient representation of $A_m$.
(Alternatively, we could use the fact that the even subgroup of 
$\langle r_j, \ldots , r_{j+m-1}\rangle$ is the alternating group of 
degree $m+1$, which is simple if $m\geq 4$; the cases $m=2,3$ are 
straightforward.) 
For $m=1$, we note that $A_1^2 = \{e\}$; 
otherwise, for $s \geq 3$, $A_1^s \simeq A_1 \simeq C_2$.

\noindent 
\textbf{(b) The group of the $m$-cube: $B_m$, for $m \geq 2$}.

We must accommodate two distinct basic systems for $B_m$.
Consider first the subdiagram
\begin{equation}\label{Bmsys1}
 \ldots \frac{}{\;\;\;}\!
\stackrel{a}{\bullet}\!
\frac{}{\;\;\;\;\;\;}\!
\stackrel{2a}{\bullet}\!
\frac{}{\;\;\;\;\;\;}\!
\stackrel{2a}{\bullet}\!
\frac{}{\;\;\;}\; 
 \ldots \; \frac{}{\;\;\;}\;
\stackrel{2a}{\bullet}\!
\frac{}{\;\;\;\;\;\;}\!
\stackrel{2a}{\bullet}\!\frac{}{\;\;\;}\;\ldots  , 
\end{equation}
on nodes $j, \ldots , j+m-1$ of $\Delta(G)$. Then 
$$ \langle r_j, \ldots , r_{j+m-1} \rangle^s \simeq B_m\;,$$
for all $s \geq 3$, \textit{and} for $s = 2$ so long as node $j$ (labelled $a$) is \oae.
If, however, $s=2$ and node $j$ is \eae, then $r_j = e$ and the given generators do not 
give a string $C$-group. Instead, from (a) we see that the subgroup collapses in rank to
a copy of $A_{m-1}^2$.  

Here, and in similar situations below,  
we obtain dual versions of these results by flipping the diagram end-for-end.
Consequently,  we may suppose that $m\geq3$ for the alternative
basic system 
\begin{equation}\label{Bmsys2} \ldots \frac{}{\;\;\;}\!
\stackrel{2a}{\bullet}\!
\frac{}{\;\;\;\;\;\;}\!
\stackrel{a}{\bullet}\!
\frac{}{\;\;\;\;\;\;}\!
\stackrel{a}{\bullet}\!
\frac{}{\;\;\;}\; 
 \ldots \; \frac{}{\;\;\;}\;
\stackrel{a}{\bullet}\!
\frac{}{\;\;\;\;\;\;}\!
\stackrel{a}{\bullet}\!\frac{}{\;\;\;}\;\ldots  , 
\end{equation}
on nodes $j, \ldots , j+m-1$.
Again we have 
$$ \langle r_j, \ldots , r_{j+m-1} \rangle^s \simeq B_m\;,$$
whenever the modulus $s \geq 3$. For $s=2$, 
the subgroup $ \langle r_j, \ldots , r_{j+m-1} \rangle^s$
is isomorphic either to $B_m$ (the group of the cube),  or to
$B_m/\{\pm e\}$ (the group of the \textit{hemi-cube}\,), as detailed 
in Table~\ref{Bmmod2}.

\begin{center}
\begin{tabular}{c|c|c|c}
\multicolumn{2}{c|}{node} & \multicolumn{2}{c}{$m \geq 3$}\\ \hline
$j$ & $ j+m -1$ & even  & odd\\ \hline\hline
\oao& \oao& $B_m$ & $B_m$ \\ \hline
\oao& \oae& $B_m/\{\pm e\}$& $B_m$\\ \hline
\oae& \oao& $B_m$& $B_m$\\ \hline
\oae& \oae& $B_m/\{\pm e\}$&$B_m/\{\pm e\}$ 
\end{tabular}
\vspace*{3mm}

\end{center}
\begin{mytab}{The group $B_m^2$ for the diagram (\ref{Bmsys2})}\label{Bmmod2}\end{mytab}
(Since $m \geq 3$, node $j+m-1$ cannot be \eae.)
Note that the bottom row covers the case
that  $G$ actually equals $B_m$,  for which there is inevitably a 
collapse when $s=2$.
A crucial step in the verification employs a small observation concerning
$B_m \simeq [4, 3, \ldots, 3] = \langle r_0, r_1, \ldots, r_{m-1} \rangle$:
if $\varphi : B_m \rightarrow H$ is a homomorphism which is $1-1$ on the subgroups
$\langle r_0, r_1 \rangle$ and 
$\langle r_1, \ldots, r_{m-1} \rangle$, then $\ker \varphi \subseteq \{\pm e \}$.
The proof follows from explicit calculation in $B_m$, taken  as the semidirect  product
$C_2^m \rtimes S_m$.
Note here that   $H$ is isomorphic to $B_m$  
if and only if 
$(r_{0}r_{1}\ldots r_{m-1})^m \neq e$. (To argue from a topological 
perspective, 
the regular $m$-polytope associated with the group 
$\langle r_j, \ldots , r_{j+m-1} \rangle^s$ must be a regular tessellation on 
an $(m-1)$-dimensional spherical space-form and hence necessarily be isomorphic 
to a regular tessellation on the $(m-1)$-sphere or real projective $(m-1)$-space 
(see \cite[6C2]{arp}). This observation also applies to the next group.)

\noindent 
\textbf{(c) The group of the $24$-cell: $F_4$}.

We must consider a  subdiagram such as 
\begin{equation}\label{F4sys}
\ldots \frac{}{\;\;\;}\!
\stackrel{a}{\bullet}\!
\frac{}{\;\;\;\;\;\;}\!
\stackrel{a}{\bullet}\!
\frac{}{\;\;\;\;\;\;}\!
\stackrel{2a}{\bullet}\!
\frac{}{\;\;\;\;\;\;}\!
\stackrel{2a}{\bullet}\!\frac{}{\;\;\;}\;\ldots  
\end{equation}
on nodes $j, \ldots , j+3$ in  $\Delta(G)$. By part (b), the natural mapping
$$\varphi : F_4 \rightarrow \langle r_j, r_{j+1} , r_{j+2}, r_{j+3} \rangle$$
is $1-1$ on subgroups $ \langle r_j, r_{j+1} , r_{j+2}\rangle$ and 
$\langle r_{j+1} , r_{j+2}, r_{j+3} \rangle$. A similar
 small observation now gives
$\ker \varphi \subseteq \{\pm e\}$.  No matter how the subdiagram is embedded
in $\Delta(G)$ we find that 
\begin{equation}\label{F4types} 
\langle r_j, r_{j+1}, r_{j+2} , r_{j+3} \rangle^s \simeq 
\left\{ \begin{array}{cl}
 F_4 \, ,& \mbox{  if } s \geq 3\, ;\\
  F_4/\{\pm e\}\, , & \mbox{ if  } s= 2\;.
\end{array}
\right.
\end{equation}

\medskip
\section{Modular polytopes  of Euclidean type}
\label{euctype}
Suppose now that  $G = [p_1, \ldots ,p_{n-1}]$ is a string Coxeter group of 
\textit{Euclidean} (or \textit{affine}) type, 
 with connected diagram (no $p_j = 2$).
Then $G$ acts as the full symmetry group of a certain regular tessellation 
$\mathcal{T} \simeq \mathcal{P}(G)$ of Euclidean  
space
 $\mathbb{A}^{n-1}$. Indeed,   $G$  
must be one of the Coxeter groups displayed
in the left column of Table~\ref{EucCoxGps}, though perhaps with generators specified in dual 
order. Note that each of these groups is  crystallographic.

A \textit{regular} $n$-\textit{toroid} $\mathcal{P}$ is the quotient of 
such a tessellation $\mathcal{T}$
 by a non-trivial normal subgroup $L$ of translations in $G$. 
Thus every toroid can be viewed 
as a finite, regular tessellation of the $(n-1)$-torus. We refer to 
\cite[1D and 6D-E]{arp} for a complete classification; briefly, for each group $G$
the distinct toroids are indexed by a \textit{type vector}
${\bf q} := (q^k, 0^{n-1-k}) = (q,\ldots,q,0,\ldots,0)$, 
where $q \geq 2$ and $k = 1, 2$ or $n-1$. 
(For $G = [3,3,4,3]$, the case $k=4$ is subsumed by the case $k=1$.)
Anyway, 
$L$ is generated (as a normal subgroup of $G$) by the translation
$$ \overline{t} := t_1^q \cdots t_k^q\;,$$
where $\{t_1, \ldots , t_{n-1}\}$ is a standard set of generators for the 
full group $T$ of translations in $G$.
The modular toroids $\mathcal{P}(G^p)$ described in 
\cite[\S 6B]{monsch1} are special instances; with one exception, we had there
${\bf q} = (p,0,\ldots, 0)$.

For completeness we also list in Table~\ref{EucCoxGps} the infinite dihedral
group $[\infty]$, which of course has rank $2$ and 
 acts on the Euclidean line 
$\mathbb{A}^1$.  The corresponding $2$-toroids are then regular polygons inscribed
in a `$1$-torus', namely,  in an ordinary circle.

Before proceeding to a classification of the groups $G^s$, 
we take a closer look at the geometric action of groups
of affine Euclidean isometries.
Suppose then that
$G = \langle r_0, \ldots , r_{n-1} \rangle $
is  of Euclidean type (here always with connected diagram). From \cite[\S 6.5]{humph}
we recall  that the invariant quadratic form
$x \cdot y$ on  real $n$-space $V$ must be positive semidefinite, so that
the \textit{radical subspace} 
$\mbox{rad}(V) = \langle c \rangle$ is $1$-dimensional. 
Since $r_j(c) = c$, for $0\leq j \leq n-1$, $G$ is in fact a subgroup of 
$\widehat{O}(V)$, the pointwise stabilizer of $\mbox{rad}(V)$ in $O(V)$. 

To actually exploit the structure of $G$ as a group of (affine)
isometries on 
Euclidean $(n-1)$-space, we pass to the contragredient representation of 
$G$ in the dual space $\check{V}$ (as in \cite[5.13]{humph}). 
Since $c$ is fixed by $G$, we see that $G$ leaves invariant any 
translate of the $(n-1)$-space
$$ U =  \{ \mu \in \check{V} \,:\, \mu(c) =0 \}\;.$$
Next,  for each $w \in V$ define $\mu_w \in \check{V}$ by 
$\mu_w(x) := w \cdot x$. The mapping
$w \mapsto \mu_w$ factors to  a linear isomorphism between 
$V/\mbox{rad}(V)$ and $U$, and  so we transfer to $U$ the 
\textit{positive definite} form induced by  $V$ on $V/\mbox{rad}(V)$. 
Now choose any $\alpha \in \check{V}$ such that $\alpha(c) = 1$, and let
$\mathbb{A}^{n-1} := U +\alpha$. Putting all this together we may now
think of $\mathbb{A}^{n-1}$ as \textit{Euclidean $(n-1)$-space},
with $U$ as its \textit{space of translations}. Indeed, each fixed $\tau \in U$
defines an isometric translation on $\mathbb{A}^{n-1}$:
$$\mu \mapsto \mu +\tau,\;\; \forall \mu \in \mathbb{A}^{n-1}\;.$$
It is easy to check that this mapping on $\mathbb{A}^{n-1}$ is induced 
by a unique 
isometry $t \in \widehat{O}(V)$, namely the \textit{transvection}
\begin{eqnarray*}
t(x) & = &x - \tau(x) c, \\
     & = &x - (x\cdot a)\, c ,
\end{eqnarray*}
where $\tau = \mu_a$ for suitable $a \in V$. (Remember  here that we employ 
the contragredient representation of $\widehat{O}(V)$ on $\check{V}$, not 
just that of $G$.)
In summary, we can therefore safely
think of translations as transvections.

In the following table we list those Euclidean Coxeter groups which are  
relevant to our analysis (see \cite[\S 6B]{monsch1}). Concerning the 
group 
$G = [4, 3^{n-3}, 4]$ (for the familar cubical tessellation of 
$\mathbb{A}^{n-1}$), we recall our convention that
$3^{n-3}$ indicates a string of $n-3 \geq 0$ consecutive $3$'s.
\begin{center}
\begin{tabular}{c|c|c|c}\label{euctypes}
 The group $G$&  $\dim(\mathbb{A}^{n-1})$& One possible diagram& 
 The corresponding vector\\
& & $\Delta(G)$ &  $c \in \mbox{rad}(V)$  \\\hline
 & &  & \\
 $[4, 3^{n-3}, 4]$ & $n-1 \geq 2$&
$ \stackrel{2}{\bullet}\!\frac{}{\;\;\;\;\;\;}\!
	\stackrel{1}{\bullet}\!\frac{}{\;\;\;\;\;\;}\!
	\stackrel{1}{\bullet} \; \cdots \;
	\stackrel{1}{\bullet}\!\frac{}{\;\;\;\;\;\;}\!
	\stackrel{1}{\bullet}\!\frac{}{\;\;\;\;\;\;}\!
	\stackrel{2}{\bullet}$
                      & $c = b_0 + 2(b_1 + \ldots + b_{n-2}) + b_{n-1}$  \\ \hline
& & &\\
$[3,3,4,3]$        & $4$& 
	$ \stackrel{1}{\bullet}\!\frac{}{\;\;\;\;\;\;}\!
	\stackrel{1}{\bullet}\!\frac{}{\;\;\;\;\;\;}\!
	\stackrel{1}{\bullet}\!\frac{}{\;\;\;\;\;\;}\!
	\stackrel{2}{\bullet}\!\frac{}{\;\;\;\;\;\;}\!
	\stackrel{2}{\bullet}$
		&  $c = b_0 + 2 b_1 + 3 b_2 + 2 b_3 + b_4$  \\ \hline
& &  &\\
$[3,6]$            &$2$ &
$ \stackrel{1}{\bullet}\!\frac{}{\;\;\;\;\;\;}\!
	\stackrel{1}{\bullet}\!\frac{}{\;\;\;\;\;\;}\!
	\stackrel{3}{\bullet}$			 & $c = b_0 + 2 b_1 + b_2$ \\ \hline
& &  &\\
$[\infty]$ & $1$ & 
$\stackrel{1}{\bullet} =\!\!=\!\!= \stackrel{1}{\bullet}$ & $c = b_0 +  b_1$ \\ 
\end{tabular}

\vspace*{4mm}
 
\end{center}

\begin{mytab}{Euclidean Coxeter Groups}\label{EucCoxGps}\end{mytab}

\noindent
An investigation of
the action of these discrete reflection groups on the Euclidean 
space $\mathbb{A}^{n-1}$ shows, in each case, that 
$G \simeq T \rtimes H$ splits as the semidirect product of 
the (normal) subgroup $T$ of translations with a certain 
(finite) \textit{point group}  $H$ (see \cite [Prop. 4.2]{humph}).
 We can and do display each  group in the table so that  
 $H = G_0 = \langle r_1, \ldots, r_{n-1} \rangle$. 

Now we are in a position to survey the 
modular reduction of the Euclidean groups in 
Table~\ref{EucCoxGps}. Again we  more generally consider Euclidean subgroups
\begin{equation}\label{embedsplit}
E = \langle r_j, \ldots, r_{j+m} \rangle \simeq T \rtimes 
\langle r_{j+1}, \ldots, r_{j+m} \rangle
\end{equation}
of  our usual group $G$; and once more we allow various
possible basic systems. Notice that we specifically assume that $E$ is embedded in
$G$ so that the point subgroup (of spherical type)
is $\langle r_{j+1}, \ldots, r_{j+m} \rangle$. Because of this,  we can use the 
splitting in (\ref{embedsplit}) to actually perform
explicit calculations, although the details 
are quite involved. We begin with

\begin{LEM}\label{toroidreduce}  
Let $G$ be a crystallographic linear Coxeter group  
with string diagram. Suppose that 
$ E = \langle r_j , \ldots, r_{j+m} \rangle$ 
is the (Euclidean) subgroup of $G$ corresponding to one of the subdiagrams
displayed in Table~\ref{Cubics} or Table~\ref{Exotics}, so that
$ E = T \rtimes H$, 
with translation group $T$ and (spherical) point group 
$H = \langle r_{j+1} , \ldots, r_{j+m} \rangle$. 
Also suppose that $s, m$, and 
the nodes $j, j+m$ are restricted in one of the various ways indicated
in the Tables, so in particular $H \simeq H^s$.
Let $\varphi : E \rightarrow E^s \subseteq G^s$ be the natural epimorphism
for  modulus $s \geq 2$. 

\mbox{\rm (a)} Then $\ker(\varphi) \subset T$.

\mbox{\rm (b)} $E^s$ is a string $C$-group, namely the automorphism
group of a regular $m$-toroid.

\mbox{\rm (c)} If $T^s$ acts faithfully on the $\mathbb{Z}_s$-submodule
spanned by $b_j , \ldots , b_{j+m}$, then 
$$T^s \cap \langle r_{j+1} , \ldots , r_{j+m}, \ldots ,  
r_{j+l}\rangle^s = \{e\}\;,$$
for any $l \geq m$.
\end{LEM}

\noindent\textbf{Proof}. As always, our calculations   may well depend
on the underlying choice of basic system
$\{ b_i\}$ for $G$, as encoded in the diagram $\Delta(G)$.  By inspection of 
the various diagrams in Tables~\ref{Cubics} and \ref{Exotics}, we 
confirm in each case that $ E = T \rtimes H$, with
$H = \langle r_{j+1} , \ldots, r_{j+m} \rangle$.
Furthermore, we also observe 
that the radical of $\sum_{k=j}^{j+m}\mathbb{R} b_k$
is spanned by an \textit{integral} vector 
$c = \sum_{k=j}^{j+m} x_{k}b_{k}$, 
in which the coefficient of $b_j$ is $x_j = 1$. 

Now for part (a)
let $ g = t h \in \ker(\varphi)$, with $t \in T, h \in H$, so that
$ t \equiv h^{-1} \pmod{s}$. For $j \leq i \leq j+m$, we have
$t(b_i) = b_i + z_i c$, with $z_i \in \mathbb{Z}$ (the coefficient of $b_{j}$ in $c$ is $1$), since $t$ is a  translation and the lattice $\sum_{k=j}^{j+m}\mathbb{Z} b_k$ is invariant under $E$;
likewise $h^{-1}(b_i) = b_i + v_i$, with $v_i  \in \sum_{k=j+1}^{j+m}\mathbb{Z} b_k$, 
since $h \in \langle r_{j+1} , \ldots, r_{j+m} \rangle$. Thus $z_i \equiv 0 \pmod{s}$,
so that $h^{-1} \equiv e \pmod{s}$. Since reduction modulo $s$
is faithful on $H$, we have $h = e$ (in characteristic $0$),
and $g = t \in T$. 

For part (b) we first of all note that 
the subgroups $H = \langle r_{j+1} , \ldots, r_{j+m} \rangle$
and $A:=\langle r_{j} , \ldots, r_{j+m-1} \rangle$  are spherical, since 
the various constraints on 
$s, m_{j, j-1}, m_{j+m, j+m+1}$ in Tables~\ref{Cubics} and \ref{Exotics}
guarantee that both subgroups are faithfully 
represented mod $s$; see Section~\ref{modsph}.
Now (b)  follows at once from (a), since $\ker(\varphi)$ is a 
normal subgroup of translations; see \cite[6D-E]{arp}. 
Here we also need to make a  
forward appeal to the computation of the type vector $\bf q$ of  
Tables~\ref{Cubics} and \ref{Exotics}, 
eliminating the possibility that the index of $\ker  
(\varphi)$ in $T$ is too small for $E^s$ to be polytopal. (We can also  
give a direct proof of the intersection property of $E^s$ using 
\cite[Prop. 2E16(a)]{arp}. 
Since the subgroups $A,H$ are both (spherical) string $C$-groups, we need only
show that $A^s \cap H^s \subseteq \langle r_{j+1} , \ldots, r_{j+m-1} \rangle^s$. 
So suppose $g \in A$ and $h \in H$ (both in  
characteristic $0$) such that $g\equiv h \bmod s$. Then $h^{-1}g =: t  
\in \ker(\varphi) \subseteq T$. Now let $\cal T$ be the regular  
tessellation in Euclidean $m$-space associated with $E$, let $o$ be  
the base vertex of $\cal T$, and let $z$ be the center of the base  
facet (tile) $F$ of $\cal T$. Then $t^{-1}(h^{-1}(z)) = g^{-1}(z) = z
$, so $t$ must be the translation by the vector $h^{-1}(z)-z$. Since  
$h^{-1}(z)$ is the center of the facet $h^{-1}(F)$ of $\cal T$ and $o
$ is a vertex of $h^{-1}(F)$, the two vertices $h^{-1}(z)$ and $z$ of   
the dual of the vertex-figure of $\cal T$ at $o$ are equivalent under  
$t$ and thus under $\ker(\varphi)$. Hence, if $t$ is non-trivial, 
then reduction modulo $s$ collapses the vertex-figure of $\cal T$ at $o$, contrary to the fact  
that $H^s$ is isomorphic to $H$. Therefore, $t$ must be trivial and
$g=h\in A\cap H = \langle r_{j+1},\ldots,r_{j+m-1}\rangle$. It follows that the modular images of 
$g$ and $h$ are in $\langle r_{j+1} , \ldots, r_{j+m-1} \rangle^s$, as required. 
Alternatively we can argue here as  
follows. The translation vectors of the conjugates of $t$ under $H$  
generate a sublattice of $\ker(\varphi)$ with very small index in $T$;  
however, our computation of the type vectors $\bf q$ has shown that  
this cannot occur.)

For part (c) we let 
$\varphi(t) = \varphi(h) \in T^s \cap \langle r_{j+1} , \ldots , r_{j+l} \rangle^s$.
Again $t(b_i) \equiv b_i\pmod{s}$ for $j \leq i \leq j+m$,
so that by hypothesis we have $t \equiv e \pmod{s}$.  \hfill $\square$ 

\noindent\textbf{Remarks}. We have seen that $H\simeq H^s$ always holds when
$s \geq 3$ and occasionally when $s=2$; under the constraints on $m$ indicated 
in Tables~4 and 5, it also holds for $s=2$.
A consequence of our calculations  is that, 
for all the cases detailed
in Tables~\ref{Cubics} and \ref{Exotics},
the semidirect splitting (\ref{embedsplit}) of 
$E =\langle r_j, \ldots, r_{j+m} \rangle$
(in characteristic $0$)  survives reduction modulo $s$.
Thus,
$E^s \simeq T^s \rtimes H^s$, although it is not necessarily the case that
$T^s \simeq \mathbb{Z}_s^m$.

Of course, taking the $r_{i}$'s in reverse order, 
we obtain a dual version of   Lemma~\ref{toroidreduce}. 
In applications, we must then take care
that the subdiagrams in Tables ~\ref{Cubics} and \ref{Exotics}, along with
the attached constraints, really have been flipped end-for-end.

Next we must deal with the specific features of each group $G$.
Guided by \cite[6D-E]{arp}, we can, with some effort, write out
explicit matrices for standard generators $t_1, \ldots, t_m$ of 
the translation subgroup $T \subset \langle r_j , \ldots, r_{j+m} \rangle$.
Such matrices incorporate the unspecified, but crucial, Cartan integers
$m_{j, j-1}$ and $m_{j+m, j+m+1}$ and furthermore vary a little 
with the choice of the underlying basic system.  But from 
Lemma~\ref{toroidreduce}(b) we know that $E^s$ is a string $C$-group. To finish
off its description, we identify the type vector
$\bf q$ by   calculating the periods of the key 
translations $t_1, t_1 t_2$ and $t_1 t_2 \ldots t_m$.  
It is convenient now to separate our results into two lots: 

\noindent
\textbf{(a) The groups $[ 4, 3^{m-2} ,4]$ ($m \geq 2$)}.

When $ \langle r_j, \ldots, r_{j+m} \rangle \simeq [ 4, 3^{m-2} ,4]$, we 
must contend
with the three distinct basic systems shown in Table~\ref{Cubics}. 
For any $s \geq 3$, we observe that   
$ \langle r_j, \ldots, r_{j+m} \rangle^s$ is the group of a suitable
cubic toroid $\{4,3^{m-2},4\}_{\bf q}$ of rank $m+1$ (on the $m$-torus), 
whose type vector ${\bf q}$ is also displayed in the Table.
The same holds for $s = 2$, so long as terminal nodes $j$ and $j+m$ are constrained as
indicated. This restriction guarantees that the facet and vertex-figure subgroups 
are spherical, with the correct rank $m$ (see Section~\ref{modsph} above).
For any other terminal node types when $s=2$, one finds that 
$ \langle r_j, \ldots, r_{j+m} \rangle^2$ either fails to have involutory 
generators (so is not a string $C$-group) or is \textit{locally projective}
rather than toroidal (see  \cite[14A]{arp} and \cite{halp}).

\begin{center}
\begin{tabular}{c|c|c|c|c}\label{cubbytypes}
  Subdiagram of $\Delta(G)$ & Modulus& Affine &  
 Constraints on& Type vector \\
  on nodes $j , \ldots , j+m$ & $s$ &dim. $ m\geq 2$&nodes $j, j+m$&  ${\bf q}$\\ \hline
     & & & \\ 
$ \frac{}{\;\;\;}\!\stackrel{2a}{\bullet}\!\frac{}{\;\;\;\;\;\;}\!
	\stackrel{a}{\bullet} \; \cdots \;
	\stackrel{a}{\bullet}\!\frac{}{\;\;\;\;\;\;}\!
	\stackrel{2a}{\bullet}\!\frac{}{\;\;\;}$
   
	     & odd $s \geq 3$ & any& --- & $(s,0, \ldots,0)$\\
     & even $s \geq 4$&  $m$ odd& at least one \oao & $(s,0, \ldots,0)$ \\
     & even $s \geq 4$&  $m$ odd &both \oae & $(\frac{s}{2}, \frac{s}{2}, \ldots, \frac{s}{2} )$\\  
     & even $s \geq 4$&  $m$ even& --- & $(\frac{s}{2}, \frac{s}{2}, \ldots, \frac{s}{2})$\\  
     &   $s =2 $&  $m$ odd &both \oao & $(2,0, \ldots, 0)$\\ \hline
     & & & &\\
 $ \frac{}{\;\;\;}\!\stackrel{a}{\bullet}\!\frac{}{\;\;\;\;\;\;}\!
	\stackrel{2a}{\bullet} \; \cdots \;
	\stackrel{2a}{\bullet}\!\frac{}{\;\;\;\;\;\;}\!
	\stackrel{a}{\bullet}\!\frac{}{\;\;\;}$   
	      & odd $s \geq 3$& any& --- &$(s,0, \ldots,0)$ \\
    & even $s \geq 4$& any & at least one \oae  &$(s,0, \ldots,0)$ \\
    & even $s \geq 4$&any & both \eae &$(\frac{s}{2},0, \ldots,0)$ \\
  &$s=2$ &any& both \oae& $(2,0, \ldots,0)$\\ \hline
    &  & & & \\
$ \frac{}{\;\;\;}\!\stackrel{4a}{\bullet}\!\frac{}{\;\;\;\;\;\;}\!
	\stackrel{2a}{\bullet} \; \cdots \;
	\stackrel{2a}{\bullet}\!\frac{}{\;\;\;\;\;\;}\!
	\stackrel{a}{\bullet}\!\frac{}{\;\;\;}$  
	   &odd $s \geq 3$ &any & ---  & $(s,0, \ldots,0)$\\
    & even $s \geq 4$ & any &   $j+m$ is  \eae & $(s,0, \ldots,0)$\\
    & even $s \geq 2$ & any  &  $j+m$ is \oae & $(s,s,0, \ldots,0)$\\ \hline
\end{tabular}

\vspace*{4mm}
 
\end{center}
 
\begin{mytab}{Groups for the cubic toroids}\label{Cubics}\end{mytab}

\noindent
\textbf{(b) The special groups $[3,3,4,3]$  ($m =4$),  $[3,6]$ ($m=2$)
 and $[\infty]$ ($m=1$)}.

Similar remarks apply to the remaining Euclidean groups 
$ \langle r_j, r_{j+1}, r_{j+2}, r_{j+3}, r_{j+4}   \rangle \simeq [3,3,4,3]$,  
$ \langle r_j, r_{j+1}, r_{j+2}  \rangle \simeq [3,6]$ or 
$ \langle r_j, r_{j+1} \rangle \simeq [\infty]$ (and their duals).
For the first two groups we may exclude the 
modulus $s=2$, for which there is a collapse in either the facet 
or vertex-figure. Our calculations   are summarized in Table~\ref{Exotics}. 
The resulting polytopes are regular toroids $\{3,3,4,3\}_{\bf q}$ of rank 
$5$ (on the $4$-torus), $\{3,6\}_{\bf q}$ of rank $3$ (on the $2$-torus), and 
regular polygons $\{q\}$ (on the $1$-torus), when ${\bf q}=(q)$ in the latter 
case.  Note
for the group $[3,6]$ that the residue of the Cartan integer $m_{j, j-1}\!\!\! \pmod{3}$ 
is a consideration (see \cite[5.6]{monsch1}).

\noindent\textbf{Remark}.  We have surveyed here the Euclidean subgroups 
$E$ of $G$. We emphasize that any reduced subgroup $E^s$ not 
explicitly covered (up to duality)
by an entry in  Table~\ref{Cubics} or Table~\ref{Exotics}
will fail in some way to be the group of a regular toroid.

\begin{center}
\begin{tabular}{c|c|c|c|c} 
  Subdiagram of $\Delta(G)$ & Modulus& Affine &  
 Constraints on& Type vector \\
  on nodes $j , \ldots , j+m$ & $s $ &dim. $m$&nodes $j, j+m$&  ${\bf q}$\\ \hline
     & & & \\ 
$ \frac{}{\;\;\;}\!\stackrel{a}{\bullet}\!\frac{}{\;\;\;\;\;}\!
	\stackrel{a}{\bullet}\!\frac{}{\;\;\;\;\;}\!
	\stackrel{a}{\bullet}\!\frac{}{\;\;\;\;\;}\!
	\stackrel{2a}{\bullet}\!\frac{}{\;\;\;\;\;}\!
	\stackrel{2a}{\bullet}\!\frac{}{\;\;\;}$
   
	     & odd $s \geq 3$ & $4$ & --- & $(s,0,0,0)$\\
     & even $s \geq 4$&  $4$ &  node $j$ is \oao & $(s,0,0,0)$ \\
     & even $s \geq 4$&  $4$ &  node $j$ is \oae & $(\frac{s}{2},\frac{s}{2},0,0)$\\ \hline 
&&&&\\
$ \frac{}{\;\;\;}\!\stackrel{2a}{\bullet}\!\frac{}{\;\;\;\;\;}\!
	\stackrel{2a}{\bullet}\!\frac{}{\;\;\;\;\;}\!
	\stackrel{2a}{\bullet}\!\frac{}{\;\;\;\;\;}\!
	\stackrel{a}{\bullet}\!\frac{}{\;\;\;\;\;}\!
	\stackrel{a}{\bullet}\!\frac{}{\;\;\;}$
   
	     & any $s \geq 3$ & $4$ & --- & $(s,0,0,0)$\\ \hline\hline 
	     &&&&\\
$ \frac{}{\;\;\;}\!\stackrel{a}{\bullet}\!\frac{}{\;\;\;\;\;}\!
	\stackrel{a}{\bullet}\!\frac{}{\;\;\;\;\;}\!
	\stackrel{3a}{\bullet}\!\frac{}{\;\;\;}$ & $s \equiv \pm 1 \pmod{3}$& 2& ---& $(s,0)$\\
	&($s>2$)&&&\\
	   &$s \equiv 0 \pmod{3}$&2& $m_{j, j-1}\equiv\pm 1\!\!\!\! \pmod{3} $ &$(s,0)$\\
	   &$s \equiv 0 \pmod{3}$&2& $m_{j, j-1}\equiv 0\!\!\! \pmod{3} $ &$(\frac{s}{3},\frac{s}{3})$\\ \hline 
	   &&&& \\
$\frac{}{\;\;\;}\!\stackrel{3a}{\bullet}\!\frac{}{\;\;\;\;\;}\!
	\stackrel{3a}{\bullet}\!\frac{}{\;\;\;\;\;}\!
\stackrel{a}{\bullet}\!\frac{}{\;\;\;}$	 &any $s \geq 3$&2& --- &  $(s,0)$\\ \hline\hline
	  &&&&\\	
$\frac{}{\;\;}\!\!\stackrel{a}{\bullet} =\!\!=\!\!= 
        \stackrel{a}{\bullet}\!\frac{}{\;\;}$ & odd $s \geq 3$ & 1 & --- & $(s)$\\
        & even $s \geq 4$ & 1 & some node \oae& $(s)$\\
        & even $s \geq 4$ & 1 & both nodes \eae & $(\frac{s}{2})$ \\
        & $s =2$ & 1 & both nodes \oae & $(2)$ \\\hline
	  &&&&\\	
$\frac{}{\;\;}\!\!\stackrel{4a}{\bullet}\!\frac{}{\;\;\;\;\;}\!
\stackrel{a}{\bullet}\!\!\frac{}{\;\;}$& odd $s \geq 3$ & 1 & --- & $(s)$\\
        & even $s \geq 4$ & 1 & node $j+1$ is \eae& $(s)$\\
        & even $s \geq 2$ & 1 & node $j+1$ is \oae& $(2s)$ \\\hline
\end{tabular}

\vspace*{4mm}
 
\end{center}

\begin{mytab}{Groups for the special toroids}\label{Exotics}\end{mytab}

\medskip
 
\section{The Quotient Criterion}\label{quotmod}

The following result is a modular variant of the \textit{quotient criterion}
in \cite[2E17]{arp}. As usual there is a dual version with subgroups
$G_{n-1}$ and $G_0$ interchanged.

\begin{THM}\label{interquotA}  
Let $G = \langle r_0, \ldots , r_{n-1} \rangle$  
be a crystallographic linear Coxeter group with string diagram, 
and suppose $G^s$ is a string $C$-group for modulus 
$s \geq 2$. 
Suppose also that $s | d$ 
and that either

\noindent {\rm (a)} $G_{n-1}$ is of spherical type
and  that $G_{n-1} \simeq G_{n-1}^s$ (so that the underlying basic system of $G$ is 
restricted as explained in \mbox{ \S\ref{modsph}} when $s=2$\,); or

\noindent {\rm (b)} $G_{n-1} = T \rtimes G_{0,n-1}$ is of Euclidean type,
with translation group $T$ and (faithfully represented) spherical point  group
$G_{0,n-1} \simeq G_{0,n-1}^s$ 
(so that $n \geq 3$ and the underlying basic system of $G$ is 
restricted as explained in \mbox{\S\ref{euctype}}\,).  Also assume in this case that
\begin{equation}\label{speceuc}
T^d \cap \langle r_1, \ldots , r_{n-1}\rangle^d = \{e\}  \;.
\end{equation}

\noindent Then  $G^d$ is a string $C$-group.
\end{THM}

\noindent\textbf{Proof}. We adapt the proof of  \cite[2E17]{arp}. 
Since $s | d$ we have natural epimorphisms $\eta : G \rightarrow G^d$ and
$\varphi : G^d \rightarrow G^s$. For clarity we avoid our customary
abuse of notation and 
take care to distinguish the standard generators $q_j := \eta(r_j)$
of $G^d$ and $s_j := \varphi(q_j)$ of $G^s$.
Since $G^s$ is a string $C$-group, 
each $s_j$ and hence each $q_j$ is an involution. By \cite[2E16(b)]{arp}, 
we need only show 
that $G_{n-1}^d$ is a string C-group and, for 
$ 1 \leq k \leq n-1$, that 
$ G_{n-1}^d \cap \langle q_k, \ldots , q_{n-1} \rangle \subseteq
\langle q_k, \ldots , q_{n-2} \rangle$. So, beginning with the latter, let 
$ g \in G_{n-1}^d \cap \langle q_k, \ldots , q_{n-1} \rangle $; then 
$\varphi(g) \in \langle s_k, \ldots , s_{n-2}  \rangle  \subseteq G_{n-1}^s$, since
$G^s$ is a string $C$-group. 
 
In the spherical case (a), $\varphi$
is 1--1 on $G_{n-1}^d$,
since $G_{n-1} \simeq G_{n-1}^s$ ($\simeq G_{n-1}^d$).
Thus $ g \in \langle q_k, \ldots , q_{n-2} \rangle$. 
 
Consider the Euclidean case (b). There exists (a unique)
$h \in \langle q_{k}, \ldots, q_{n-2} \rangle$ with $\varphi(h) = \varphi(g)$.
Applying Lemma~\ref{toroidreduce} to 
$\varphi \circ \eta$ (restricted to $G_{n-1}$),  we have
$g = t h$ for some translation $t \in T^d$. By (\ref{speceuc}) we get $t = e$,
so that $g \in \langle q_k, \ldots , q_{n-2} \rangle$.

Finally, $G_{n-1}^d$ is a string C-group in each case. This follows from applying
our considerations in 
Sections \ref{modsph} and  \ref{euctype} to $G_{n-1}$, since switching 
from $s$ to a multiple $d$
merely eases any constraints which could prevent
$G_{n-1}^d$ from being a string $C$-group. \hfill $\square$

\noindent\textbf{Example and Remarks}. In general, some condition like
(\ref{speceuc}) is necessary. Consider, for instance, the  
diagrams
\begin{equation}\label{rk3egs}
\begin{array}{ccc}
\stackrel{a}{\bullet}\!\frac{}{\;\;\;\;\;}\!
	\stackrel{1}{\bullet} =\!\!=\!\!= \stackrel{1}{\bullet} 
\;&\;\mbox{ and } \;&\;
\stackrel{1}{\bullet}\!\frac{}{\;\;\;\;\;}\!
	\stackrel{a}{\bullet} =\!\!=\!\!= \stackrel{a}{\bullet}\;. \\
	\end{array}
\end{equation}
For $a \in \{1,2,3\}$, the  
corresponding groups of rank $3$  
reduce to  string $C$-groups for any modulus $d>2$. 
In   the left diagram  we can even take $a=4$ and so obtain a polyhedron
of type $\{d,d\}$, for odd $d \geq 3$, or type
$\{ d , \frac{d}{2} \}$, for even $d \geq 4$. However, taking $a = 4$ in  
the right diagram, we find that the intersection condition 
fails precisely when the modulus
$d = 2s$, with $s$ odd: for then $t = (r_0 r_1)^s = (r_1 r_2)^s \neq e
 \pmod{d}$;
and  $t \in T^d  \cap \langle r_1, r_2 \rangle^d$  directly
contradicts (\ref{speceuc}). We shall see that the fault lies in 
the embedding of the subdiagrams for facet and vertex-figure.

To explain what is going on  we   use  Lemma~\ref{toroidreduce}(c)  
(with $s=d, j=0, m=n-2,$ $l=m+1$). Thus 
(\ref{speceuc}) is fulfilled
whenever $T^d$ acts faithfully on the $\mathbb{Z}_d$-submodule 
$V_{n-1}$.
This holds, for example, when dropping  
node $n-1$ has no effect on  the  
embedding constraints for $G_{n-1}$, as  described in the Tables.
To see this, note  that $r_i$ 
induces a mapping $\tilde{r}_i$ on $V_{n-1}$, for $0 \leq i \leq n-2$.
Clearly, $K^d:= \langle \tilde{r}_0, \ldots , \tilde{r}_{n-2}\rangle$
is just the (toroidal) group corresponding
to the the subdiagram of $\Delta(G)$ obtained by deleting node $n-1$.
If, as we suppose,  this deletion has no effect  on the constraints
on node $n-2$, it must be that $G_{n-1}^d$ and $K^d$
have the same type vector $\bf q$, as given in the Tables. Since 
the corresponding spherical point groups are isomorphic, it follows that
$G_{n-1}^d \simeq K^d$ and that $T^d$ acts faithfully on $V_{n-1}$.
Thus $G^d$ is a string $C$-group.
In particular, we now see that 
(\ref{speceuc}) is redundant whenever $d$ is odd and 
in several other instances. This leads to an important simplification:
for $d$ odd we need only check   that $G^s$ is a string $C$-group 
for some \textit{odd prime} divisor $s = p$. Occasionally, the modulus
$s=4$ is another keystone.

\section{Locally toroidal polytopes}
\label{loctor}

In this Section, we consider \textit{locally toroidal}  regular polytopes, 
that is polytopes
of rank $n \geq 4$
whose facets and vertex-figures are globally spherical or toroidal, 
as described above,
with at least one kind toroidal. 
The $n$-polytopes of this kind have not yet been fully classified, although 
quite a lot is known (see
\cite[Chs. 10-12]{arp}). 

As usual, we begin with a crystallographic linear Coxeter group 
$G = \langle r_0, \ldots , r_{n-1} \rangle$, but 
immediately  
discard degenerate cases in which the underlying
diagram $\Delta(G)$ is  disconnected. (In such cases
$G$ is reducible; 
and $\mathcal{P}(G)$ and its quotients have the
sort of `flatness' described in
\cite[4E]{arp}.)

In \cite{monsch2} we  discussed all locally toroidal $4$-polytopes 
$\mathcal{P}(G^p)$ which arise from our construction with prime modulus $p$. 
Since our methods for general moduli $s$
add little to the discussion of such polytopes
in \cite[Chs. 10--11]{arp} and \cite{monsch2}, we   examine here just 
one group of rank $4$, namely $G = [3,6,3]$, with diagram
$$\stackrel{3}{\bullet}\!\frac{}{\;\;\;\;\;\;}\!
	\stackrel{3}{\bullet}\!\frac{}{\;\;\;\;\;\;}\!
        \stackrel{1}{\bullet}\!\frac{}{\;\;\;\;\;\;}\!
	\stackrel{1}{\bullet}
\;.$$
When $s=4$ we find that $G^4$ has order $7680$ and is the automorphism
group of a locally toroidal $4$-polytope in the class
$\langle\,   \{3,6\}_{(4,0)} , \{6,3\}_{(4,0)} \, \rangle$. Next we note in
Table~\ref{Exotics} that there are no embedding contraints on node $2$. We
conclude from Theorem~\ref{interquotA}(b) (and the subsequent remarks) and
from \cite[p. 345]{monsch2} that $G^d$ is a string $C$-group
whenever the modulus $d$ is divisible by either $4$ or an odd prime, that is,
whenever $d \geq 3$.  The polytope $\mathcal{P}(G^d)$
is in the class 
$\langle\,   \{3,6\}_{\bf q} , \{6,3\}_{\bf r} \, \rangle$, where
always ${\bf q} = (d,0)$, but ${\bf r} = (d,0)$  when $ 3 \nmid d$ and 
${\bf r} = (\frac{d}{3},\frac{d}{3})$ when $3 \mid d$.
This construction  complements the approach
in \cite[11E]{arp}.

Turning to higher rank 
$n>4$, we observe that any spherical
facet, or vertex-figure, must  be of type $\{3^{n-2}\}$, $\{4, 3^{n-3}\}$, 
$\{ 3^{n-3},4\}$
or $\{3,4,3\}$ ($n=5$ only). Likewise,  the required Euclidean section must 
have type 
$\{4, 3^{n-4}, 4\}$ or when $n=6$, $\{3,3,4,3\}$ or $\{3,4,3,3\}$.
As described in \cite[Lemma 10A1]{arp}, these constraints severely limit 
the possibilities:
in rank $5$, we have just  $G = [4,3,4,3]$ acting on hyperbolic space 
$\mathbb{H}^4$; and
in rank $6$ we have $G = [4,3,3,4,3], [3,4,3,3,3]$ or $[3,3,4,3,3]$,
all acting on $\mathbb{H}^5$.
Thus we may complete our discussion by examining the modular polytopes which
result from these groups in ranks $5$ and  $6$.

\medskip

\subsection{Rank $\mathbf{5}$: the group $G = [4,3,4,3]$ }

Here we must contend with the four distinct basic systems
encoded in the diagrams

\begin{equation}\label{5loctordgs}
\begin{array}{ccc}
\stackrel{1}{\bullet}\!\frac{}{\;\;\;\;\;\;}\!
	\stackrel{2}{\bullet}\!\frac{}{\;\;\;\;\;\;}\!
        \stackrel{2}{\bullet}\!\frac{}{\;\;\;\;\;\;}\!
	\stackrel{4}{\bullet}\!\frac{}{\;\;\;\;\;\;}\!
	\stackrel{4}{\bullet}& &
	\stackrel{1}{\bullet}\!\frac{}{\;\;\;\;\;\;}\!
	\stackrel{2}{\bullet}\!\frac{}{\;\;\;\;\;\;}\!
        \stackrel{2}{\bullet}\!\frac{}{\;\;\;\;\;\;}\!
	\stackrel{1}{\bullet}\!\frac{}{\;\;\;\;\;\;}\!
	\stackrel{1}{\bullet}\\
(a)& &(b)\\ 
\stackrel{2}{\bullet}\!\frac{}{\;\;\;\;\;\;}\!
	\stackrel{1}{\bullet}\!\frac{}{\;\;\;\;\;\;}\!
        \stackrel{1}{\bullet}\!\frac{}{\;\;\;\;\;\;}\!
	\stackrel{2}{\bullet}\!\frac{}{\;\;\;\;\;\;}\!
	\stackrel{2}{\bullet}& &
	\stackrel{4}{\bullet}\!\frac{}{\;\;\;\;\;\;}\!
	\stackrel{2}{\bullet}\!\frac{}{\;\;\;\;\;\;}\!
        \stackrel{2}{\bullet}\!\frac{}{\;\;\;\;\;\;}\!
	\stackrel{1}{\bullet}\!\frac{}{\;\;\;\;\;\;}\!
	\stackrel{1}{\bullet}\\
(c) & & (d)
\end{array}
\end{equation}
When the modulus is an odd prime $p$, the four corresponding 
finite groups $G^p$
 are isomorphic string $C$-groups; and we 
recall from \cite[\S 4.1]{monsch3} that
\begin{equation}\label{loctor5}
G^p =
\left\{ \begin{array}{ll}
O_1(5,p,0)\;,  & \mbox{ if }  p \equiv \pm 1  \pmod{8}\\
O(5,p,0)\;,  & \mbox{ if }  p \equiv \pm 3 \pmod{8}
\end{array}\right.
\end{equation}
Note that $O_1(5,p,0)$ has order 
$p^4 (p^4 -1) (p^2 - 1) $ and index two in $O(5,p,0)$ 
(see \cite[pp. 300-301]{monsch1}).
The facets of the corresponding 
regular $4$-polytope $\mathcal{P}(G^p)$
 are toroids  $\{4,3,4\}_{(p,0,0)}$, which one 
could construct by identifying opposite square faces of a 
$p \times p \times p$ cube \cite[6.4]{monsch1}. Of course, 
the vertex-figures are copies of the $24$-cell $\{3,4,3\}$.

Next, for modulus $s=4$, we may check directly on GAP
that $G^4$ is a string $C$-group for each of
the basic systems in (\ref{5loctordgs}). Diagrams (a), (b), (c)
give polytopes of type $\{\,\{4,3,4\}_{(4,0,0)}\, , \, \{3,4,3\}\,\}$,
whose respective automorphism groups have orders
$g = 2^{16} \cdot 3^2$, $g$, and $4g$. On the other hand, diagram
(d) gives a polytope of type
$\{\,\{4,3,4\}_{(4,4,0)}\, , \, \{3,4,3\}\,\}$ whose group has order $16g$.
By \cite[12B1]{arp}, none of these polytopes can be universal for their type.
However, with different generators, the third group,
of order $4g = 2\,359\,296$, is the automorphism group for 
the universal polytope
of type $\{\,\{4,3,4\}_{(2,2,2)}\, , \, \{3,4,3\}\,\}$ and
hence is known to be isomorphic to 
$(\mathbb{Z}_{2}^{6} \rtimes \mathbb{Z}_{2}^{5}) \rtimes F_4$ 
(see \cite[Thm. 8F19 and Table 12B1]{arp}).

Now consider \textit{any} modulus $d>2$, which again is divisible either
by an odd prime $s$ or by $s =4$. We immediately conclude from
Theorem~\ref{interquotA}(a), in its dual form, that
$G^d$ is a string $C$-group for each diagram 
in (\ref{5loctordgs}) and for each modulus $d>2$.

If $d$ is \textit{odd},
it is easy to check  that the four diagrams deliver isomorphic groups.
Indeed, a change from any one of the four
basic systems  to another is accomplished by rescaling various $b_j$'s by powers 
of $2$ (see \cite[p. 305]{monsch1}). Since $2$ is invertible modulo $d$, the 
corresponding linear groups are conjugate in $GL_5(\mathbb{Z}_d)$; and,
crucially, such isomorphisms pair off the specified generating reflections.
 Consulting Table~\ref{Cubics} (with $s$ replaced by $d$), we conclude 
that the resulting non-universal
polytope has type
\begin{equation}\label{5types}
\{\,\{4,3,4\}_{(d,0,0)}\, , \, \{3,4,3\}\,\}\;.
\end{equation} 

For $d$ \textit{even}, we have already observed that a change in basic system may well
alter the corresponding group and polytope. Referring again to Table~\ref{Cubics},
we do find that diagrams (a), (b), (c) in (\ref{5loctordgs})   provide
polytopes of the type displayed in (\ref{5types}),   now with $d$ even.
However,    diagram (\ref{5loctordgs})(d)  gives a polytope of type
$\{\,\{4,3,4\}_{(d,d,0)}\, , \, \{3,4,3\}\,\}$.

Of course, in all the above cases, we just as easily 
obtain the dual polytope
of type $\{3,4,3,4\}$ by flipping a diagram end-for-end.

The universal locally toroidal polytopes of rank $5$ are  
described in \cite[12B]{arp}. 
There are just three finite instances, whose  facets are toroids with type vector 
$(2,0,0)$, $(2,2,0)$ or $(2,2,2)$. Unfortunately, we cannot get any of these
by our construction, since for $s=2$  we always have by
(\ref{F4types}) that the
24-cell collapses to its central quotient, the 
`hemi-24-cell' $\{3,4,3\}_6$.  On the other hand, for $d>2$ our
construction gives finite polytopes of the type indicated; in contrast,
the methods in \cite[p. 452]{arp} are non-constructive and appeal
to the residual finiteness of certain groups to establish
the existence of such polytopes. 

Finally, in this subsection,   
it is of some interest to further investigate the case $s=2$.
We may discard diagrams (a) and (b), in which $r_0 = e \pmod{2}$. 
However, diagram (c) does give a string $C$-group $G^2$ of order
2304, for the universal polytope
$$
\{\, \mathcal{K}\, , \, \{3,4,3\}_6\,\}\; ,
$$
where $\mathcal{K}:=\{\,\{4,3\}_3\, , \,\{3,4\} \,\}$, so that
$3$-faces and vertex figures
are of projective type. Diagram (d) likewise gives
a group $G^2$ of order 9216; and the corresponding polytope 
is doubly covered by the universal polytope of type
$$\{\,\{4,3,4\}_{(2,2,0)}\, , \, \{3,4,3\}_6\,\}\;,$$
whose group is $\mathbb{Z}_{2}^{5} \rtimes (F_{4}/\{\pm e\})$ 
(see \cite[Thm. 8F21]{arp}).

\vspace*{5mm}
\subsection{Rank $\mathbf{6}$: the groups $[3,4,3,3,3]$, $[3,3,4,3,3]$ and 
$[4,3,3,4,3]$}

In rank $6$ we must consider three closely related groups, beginning with
$$ G  = \langle r_0, r_1, r_2, r_3, r_4, r_5\rangle 
\simeq [3,4,3,3,3]\;.$$
A 
basic system (of roots) 
for $G$ is described by one of the following diagrams:

\begin{equation}\label{loctor6A}
\begin{array}{ccc}
\stackrel{1}{\bullet}\!\frac{}{\;\;\;\;\;\;}\!
	\stackrel{1}{\bullet}\!\frac{}{\;\;\;\;\;\;}\!
	\stackrel{2}{\bullet}\!\frac{}{\;\;\;\;\;\;}\!  
	\stackrel{2}{\bullet}\!\frac{}{\;\;\;\;\;\;}\!
	\stackrel{2}{\bullet}\!\frac{}{\;\;\;\;\;\;}\!
	\stackrel{2}{\bullet} & & 
\stackrel{2}{\bullet}\!\frac{}{\;\;\;\;\;\;}\!
	\stackrel{2}{\bullet}\!\frac{}{\;\;\;\;\;\;}\!
	\stackrel{1}{\bullet}\!\frac{}{\;\;\;\;\;\;}\!  
	\stackrel{1}{\bullet}\!\frac{}{\;\;\;\;\;\;}\!
	\stackrel{1}{\bullet}\!\frac{}{\;\;\;\;\;\;}\!
	\stackrel{1}{\bullet}\;\;.\\ 
(a) & & (b)	
\end{array} 
\end{equation}
Next we turn to
the subgroup $H = \langle s_0, \ldots , s_5 \rangle$ 
generated by the reflections 
\begin{equation}\label{Hgens}
(s_0, s_1, s_2, s_3, s_4, s_5) := (r_1, r_0, r_2 r_1 r_2, r_3, r_4, r_5)\; ,
\end{equation} 
which has index $5$ in $G$ and
is isomorphic to $[3,3,4,3,3]$. 
Starting with the diagram (\ref{loctor6A})(b), we find that 
the basic system of roots attached to the $s_j$'s is now encoded
in the diagram
\begin{equation}\label{loctor6B}
\stackrel{2}{\bullet}\!\frac{}{\;\;\;\;\;\;}\!
	\stackrel{2}{\bullet}\!\frac{}{\;\;\;\;\;\;}\!
	\stackrel{2}{\bullet}\!\frac{}{\;\;\;\;\;\;}\!  
	\stackrel{1}{\bullet}\!\frac{}{\;\;\;\;\;\;}\!
	\stackrel{1}{\bullet}\!\frac{}{\;\;\;\;\;\;}\!
	\stackrel{1}{\bullet}\;\;. 
\end{equation}
(Diagram 
(\ref{loctor6A})(a) merely leads, in dual fashion, to  (\ref{loctor6B})
flipped end-for-end. This is the only other diagram admitted by $H$.)

The final subgroup
$K = \langle t_0, \ldots , t_5 \rangle$ generated by
\begin{equation}\label{Kgens}
(t_0, t_1, t_2, t_3, t_4, t_5) := 
(r_2, r_1, r_0, r_3 r_2 r_1 r_2 r_3, r_4, r_5)\end{equation}
has index $10$ in $G$ and  is isomorphic to $[4,3,3,4,3]$. Now
diagrams (\ref{loctor6A})(a),(b) lead to diagrams 
(\ref{loctor6C})(a),(b) below:
\begin{equation}\label{loctor6C}
\begin{array}{ccc}
\stackrel{2}{\bullet}\!\frac{}{\;\;\;\;\;\;}\!
	\stackrel{1}{\bullet}\!\frac{}{\;\;\;\;\;\;}\!
	\stackrel{1}{\bullet}\!\frac{}{\;\;\;\;\;\;}\!  
	\stackrel{1}{\bullet}\!\frac{}{\;\;\;\;\;\;}\!
	\stackrel{2}{\bullet}\!\frac{}{\;\;\;\;\;\;}\!
	\stackrel{2}{\bullet} & &
\stackrel{1}{\bullet}\!\frac{}{\;\;\;\;\;\;}\!
	\stackrel{2}{\bullet}\!\frac{}{\;\;\;\;\;\;}\!
	\stackrel{2}{\bullet}\!\frac{}{\;\;\;\;\;\;}\!  
	\stackrel{2}{\bullet}\!\frac{}{\;\;\;\;\;\;}\!
	\stackrel{1}{\bullet}\!\frac{}{\;\;\;\;\;\;}\!
	\stackrel{1}{\bullet}\\
	(a)& & (b)\\
\stackrel{4}{\bullet}\!\frac{}{\;\;\;\;\;\;}\!
	\stackrel{2}{\bullet}\!\frac{}{\;\;\;\;\;\;}\!
	\stackrel{2}{\bullet}\!\frac{}{\;\;\;\;\;\;}\!  
	\stackrel{2}{\bullet}\!\frac{}{\;\;\;\;\;\;}\!
	\stackrel{1}{\bullet}\!\frac{}{\;\;\;\;\;\;}\!
	\stackrel{1}{\bullet} & &
\stackrel{1}{\bullet}\!\frac{}{\;\;\;\;\;\;}\!
	\stackrel{2}{\bullet}\!\frac{}{\;\;\;\;\;\;}\!
	\stackrel{2}{\bullet}\!\frac{}{\;\;\;\;\;\;}\!  
	\stackrel{2}{\bullet}\!\frac{}{\;\;\;\;\;\;}\!
	\stackrel{4}{\bullet}\!\frac{}{\;\;\;\;\;\;}\!
	\stackrel{4}{\bullet}\\
	(c)& & (d)\\	
\end{array}
\end{equation}
The group $K$ admits  the two other basic systems
shown  in (\ref{loctor6C})(c),(d).
(See \cite[12A2]{arp}.  Each group described above
acts on $\mathbb{H}^5$ with a 
simplicial fundamental domain of finite volume. In \cite{sqf}, 
these indices were computed by  
dissecting   a simplex for $H$ (or $K$)  into copies of the simplex for $G$.)

In \cite[\S 4.2]{monsch3} we showed that $G^p, H^p, K^p$ are 
string $C$-groups for any odd prime modulus $p$. In fact,
 all three are isomorphic to
\begin{equation}\label{loctor6gp}
\left\{ \begin{array}{ll}
O_1(6,p,+1)\;,  & \mbox{ if }  p \equiv \pm 1  \pmod{8}\\
O(6,p,+1)\;,  & \mbox{ if }  p \equiv \pm 3 \pmod{8}
\end{array}\right.
\end{equation}
Of course, we require different generators in the three cases,
as indicated  in (\ref{Hgens}) and (\ref{Kgens}). 
Thus, the indices $5$ and $10$ in characteristic $0$
collapse to $1$ under reduction 
$\bmod \, p$.
For any prime $p \geq 3$, $O_1(6,p,+1)$ has order 
$p^6 (p^4 -1)(p^3 - 1)(p^2 - 1) $ and index 
two in $O(6,p,+1)$ (see \cite[pp. 300-301]{monsch1}).

Now suppose that the modulus is any \textit{odd} integer $d \geq 3$. Just 
as in the previous
 subsection, the two diagrams in (\ref{loctor6A})
give isomorphic  groups, as do the four diagrams in (\ref{loctor6C}).
Furthermore, by the remarks following Theorem~\ref{interquotA} 
we see that $G^d$, $H^d$ and $K^d$ are then string $C$-groups. In each case, 
the type vector for a  toroidal section is ${\bf q} = (d,0,0,0)$.

The situation for \textit{even} moduli is more complicated. Once more, we may 
discard the modulus $d=2$, which invariably causes a collapse to the 
hemi-24-cell in any section of type $\{3,4,3\}$.  Let us consider the three 
groups in turn.

\noindent\textbf{The Polytopes $\mathcal{P} = \mathcal{P}(G^d)$}.

Using GAP, we find that 
$G^4$ is a string $C$-group of order $2^{26} \cdot 3^2 \cdot 5$ for either
diagram in (\ref{loctor6A}). It follows from 
Theorem~\ref{interquotA}(a) in its dual form that 
$G^d$ is a string $C$-group for any modulus $d>2$. 
From either diagram  in (\ref{loctor6A}) we obtain
a locally toroidal polytope
in the class
$$ \langle \{3,4,3,3 \}_{(d,0,0,0)} \, , \, \{ 4,3,3,3 \}
\rangle\;.
$$
We note that the toroidal facets of $\mathcal{P}(G^d)$ each have
$3 d^4$ vertices \cite[Table 6E1]{arp}; and, of course,  
the vertex-figures are $5$-cubes
$\{4,3,3,3\}$.  
Although the two admissible diagrams do yield string $C$-groups,
we have no general proof that these  groups are isomorphic
when $d$ is \textit{even}, though this is true for $d = 4$.

The following  theorem establishes \cite[Conjecture 12C2]{arp} concerning the 
existence of locally toroidal regular $6$-polytopes of 
type $\{3,4,3,3,3\}$.

\begin{THM}
The universal regular $6$-polytopes 
$\{\{3,4,3,3\}_{(d,0,0,0)},\{4,3,3,3\}\}$\\ and
$\{\{3,4,3,3\}_{(d,d,0,0)},\{4,3,3,3\}\}$ exist for all $d\geq 2$.
\end{THM}

\noindent\textbf{Proof}.
First note that the case $d=2$ was settled 
in \cite[pp.460-461]{arp}. So let $d>2$. We now appeal to our earlier remark 
that a non-empty class of regular polytopes 
contains a (unique) universal member (see \cite[4A2]{arp}). 
Thus, the existence of a universal polytope of the 
first kind (type vector ${\bf q} = (d,0,0,0)$)
follows directly from our construction of 
a member of its class, namely ${\cal P}(G^d)$. For the 
existence of the universal polytopes of the second kind 
(type vector ${\bf q} = (d,d,0,0)$)
we  refer to the discussion in \cite[pp.460-462]{arp}, where 
it was shown that the existence of the universal polytopes
of the second kind is implied by   existence of   
universal polytopes of the first kind. (In fact, some 
of the arguments provided there can now be simplified 
using properties of $G^d$.)
\hfill $\square$
 
The full classification of the finite universal polytopes of each 
kind is still open, but three of these  are known to be finite, 
including 
\[ \{\{3,4,3,3\}_{(3,0,0,0)},\{4,3,3,3\}\},\] 
with automorphism group
$\mathbb{Z}_3 \rtimes O(6,3,+1)\; (= \mathbb{Z}_3 \rtimes G^3)$.
See \cite[\S 4.2]{monsch3}. 

\medskip

\noindent\textbf{The Polytopes $\mathcal{P} = \mathcal{P}(H^d)$}.

We have already indicated that for $d$ odd  
 the polytope $\mathcal{P}(H^d)$  lies in the class
$$\displaystyle{\langle
\{3,3,4,3 \}_{(d,0,0,0)} \, , \, \{ 3,4,3,3\}_{(d,0,0,0)}
\rangle}\;.$$ 
In fact, $\mathcal{P}(H^d)$ admits an order reversing bijection and so 
is \textit{self-dual}.

The modulus $p =3$ is of particular interest.
In \cite[\S 4.2]{monsch3} we  gave a new construction for the corresponding
(finite!) self-dual universal polytope
$$ \mathcal{U}_{H^3} := \{ \,\{3,3,4,3\}_{(3,0,0,0)} \, , 
\, \{3,4,3,3\}_{(3,0,0,0)}\,  \}\;.$$
Indeed,
$\Gamma(\mathcal{U}_{H^3}) \simeq (\mathbb{Z}_3 \oplus \mathbb{Z}_3)
\rtimes H^3$ under a non-trivial action of $H^3$ on the abelian factor.
Thus $\mathcal{U}_{H^3}$ is a $9$-fold cover of 
$\mathcal{P}(H^3)$ (\cite[Table 12D1]{arp}); and trapped between we find
a twin pair $\mathcal{Q}, \mathcal{Q}^*$
of \textit{non-self-dual} polytopes, with the same toroidal facets
and vertex-figures:

$$
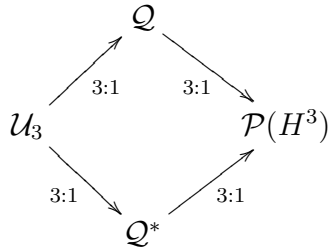
\xymatrix{ & {\mathcal{Q}}\ar[dr]_{3 : 1}&\\
          {\mathcal{U}_{3}}\ar[ur]_{3 : 1} \ar[dr]_{3 : 1} & 
          & {\mathcal{P}(H^3)}\\
           & {\mathcal{Q}^*} \ar[ur]_{3 : 1}& } 
$$

Turning to even moduli, we again 
 find that $H^4$ is a string $C$-group (of index $5$ in $G^4$); 
 and we note that there
are no embedding constraints on node $4$ (look at the second
diagram in Table~\ref{Exotics}). Thus, 
by the discussion
following Theorem~\ref{interquotA}, we conclude that $H^d$ is a string $C$-group
for all $d>2$. When $d$ is even, the corresponding polytope is in the class
$$ 
\langle \{3,3,4,3 \}_{(d,0,0,0)} \, , 
\, \{ 3,4,3,3 \}_{(\frac{d}{2},\frac{d}{2},0,0)}  
\rangle\;,
$$ 
and hence is certainly not self-dual. 

Notice that the type 
vectors for the facets and vertex-figures of the polytopes $\mathcal{P}(H^d)$ 
are related in that they involve the same parameter $d$. Thus we cannot 
expect our methods to completely
settle Conjecture 12D3 of \cite{arp} concerning
the existence of locally toroidal regular $6$-polytopes of types 
$\{3,3,4,3,3\}$, for which the parameters for the facets and 
vertex-figures may vary independently. The same remark applies 
to the polytopes $\mathcal{P}(K^d)$ studied next, and Conjecture 12E3 of \cite{arp} for the corresponding type $\{4,3,3,4,3\}$.

\noindent\textbf{The Polytopes $\mathcal{P} = \mathcal{P}(K^d)$}.

For odd $d \geq 3$ the four diagrams in (\ref{loctor6C})
give isomorphic polytopes in the class

$$ 
\langle \{4,3,3,4 \}_{(d,0,0,0)} \, , 
\, \{ 3,3,4,3 \}_{(d,0,0,0)}
\rangle\;.
$$
Here the facets are cubical toroids; facets and vertex-figures
each have  $d^4$ vertices.

Suppose then that $d \geq 4$  is even. A calculation with GAP 
reveals the at first surprising result that the intersection condition
(\ref{geninter}) fails  for diagrams (\ref{loctor6C})(b)(d),
at least when $d = 4,6$. Noting that dropping the last node in
each case
alters the constraints on node $4$, we therefore abandon these diagrams.

For diagram (\ref{loctor6C})(a) we easily verify that $K^4$ is a 
string $C$-group (of index $10$ in $G^4$). Note that there are no 
embedding constraints 
on node $4$; see the first diagram in Table~\ref{Cubics}, with $m=4$
and $s = d$ even. From Theorem~\ref{interquotA}, we thus obtain a polytope
in the class
$$ 
\langle \, \{ 4,3,3,4 \}_{(\frac{d}{2},\frac{d}{2},\frac{d}{2},\frac{d}{2} )}
\{3,3,4,3 \}_{(\frac{d}{2},\frac{d}{2},0,0)} \, \rangle\;,\mbox{ ( even $d\geq 4$ )}.
$$
Here the facets have $d^4/2$ vertices; and each  vertex-figure 
has $d^4/4$ vertices.

The analysis for diagram  (\ref{loctor6C})(c) 
is similar, although the particular location of the subgroup $[3,4,3]$ 
prevents an automatic verification of 
condition (\ref{speceuc}). Nevertheless, by brute-force 
calculation, we find that (\ref{speceuc}) 
holds for any modulus $d \geq 2$.
On the other hand, for $d = 4$ with this basic system, we can independently
check on GAP that $K^4$  is indeed a string $C$-group, with (unexpected) order 
$2^{29} \cdot 3^2$. It follows from 
Theorem~\ref{interquotA}(b) 
that $K^d$ is a  string $C$-group for any modulus $d \geq 3$.
In particular, when $d\geq 4$ is even we obtain a polytope in the class
$$ 
\langle \, \{ 4,3,3,4 \}_{(d,d,0,0)}
\{3,3,4,3 \}_{(d,0,0,0)} \, \rangle\;  .
$$
Here the facets each have $2d^4$ vertices.


\noindent
{\bf Acknowledgement.\/}
We wish to thank an anonymous referee for the careful reading of the manuscript and several helpful comments on it.

\end{document}